\documentclass[12pt,leqno]{article}
\usepackage{amsmath, amsthm, amsfonts, amssymb, color}
\usepackage{mathrsfs}
\theoremstyle{definition}

\newcommand{\scr}[1]{\mathscr #1}
\definecolor{wco}{rgb}{0.5,0.2,0.3}

\numberwithin{equation}{section} \theoremstyle{remark}

\newcommand{\ua}{\uparrow}

\title{{\bf Harnack Inequalities for Functional SDEs with Multiplicative Noise and Applications}\footnote{Supported in part by WIMCS and SRFDP}
}
\author{
{\bf Feng-Yu Wang$^{a),b)}$ and Chenggui Yuan$^{b)}$}\\
\footnotesize{$^{a)}$School of Mathematical Sciences,
Beijing Normal
University, Beijing 100875, China}\\
 \footnotesize{$^{b)}$Department of Mathematics,
Swansea University, Singleton Park, SA2 8PP, UK}\\ \footnotesize{wangfy@bnu.edu.cn, F.-Y.Wang@swansea.ac.uk, C.Yuan@swansea.ac.uk}
}
\begin{document}
\def\R{\mathbb R}  \def\ff{\frac} \def\ss{\sqrt} \def\BB{\mathbb
B}
\def\N{\mathbb N} \def\kk{\kappa} \def\m{{\bf m}}
\def\dd{\delta} \def\DD{\Delta} \def\vv{\varepsilon} \def\rr{\rho}
\def\<{\langle} \def\>{\rangle} \def\GG{\Gamma} \def\gg{\gamma}
  \def\nn{\nabla} \def\pp{\partial} \def\tt{\tilde}
\def\d{\text{\rm{d}}} \def\bb{\beta} \def\aa{\alpha} \def\D{\scr D}
\def\E{\mathbb E} \def\si{\sigma} \def\ess{\text{\rm{ess}}}
\def\beg{\begin} \def\beq{\begin{equation}}  \def\F{\scr F}
\def\Ric{\text{\rm{Ric}}} \def\Hess{\text{\rm{Hess}}}\def\B{\scr B}
\def\e{\text{\rm{e}}} \def\ua{\underline a} \def\OO{\Omega} \def\sE{\scr E}
\def\oo{\omega}     \def\tt{\tilde} \def\Ric{\text{\rm{Ric}}}
\def\cut{\text{\rm{cut}}} \def\P{\mathbb P} \def\ifn{I_n(f^{\bigotimes n})}
\def\C{\scr C}      \def\aaa{\mathbf{r}}     \def\r{r}
\def\gap{\text{\rm{gap}}} \def\prr{\pi_{{\bf m},\varrho}}  \def\r{\mathbf r}
\def\Z{\mathbb Z} \def\vrr{\varrho} \def\ll{\lambda}
\def\L{\scr L}\def\Tt{\tt} \def\TT{\tt}\def\II{\mathbb I}
\def\i{{\rm i}}\def\Sect{{\rm Sect}} \def\Q{\mathbb Q} \def\supp{{\rm supp}}

\maketitle
\begin{abstract}   By constructing  a new coupling, the log-Harnack inequality is established for the functional solution of a delay stochastic differential equation with multiplicative noise. As applications, the strong Feller property and heat kernel estimates w.r.t. quasi-invariant probability measures are derived for the associated transition semigroup of the solution. The dimension-free Harnack inequality in the sense of \cite{W97} is also investigated.
\end{abstract} \noindent
 AMS subject Classification:\  60H10, 47G20.   \\
\noindent
 Keywords:   Harnack inequality, functional solution, delay SDE, strong Feller property, heat kernel.
 \vskip 2cm

\section{Introduction}
 The dimension-free Harnack inequality introduced in \cite{W97} has become a useful tool in the study of diffusion semigroups, in particular, for the uniform integrability, contractivity properties, and estimates on heat kernels, see e.g. \cite{ATW, AZ, DRW09, K, L, LW, W99, W07, WX10,Z} and references within.  Recently, by using coupling arguments,  the dimension-free Harnack inequality has been established in \cite{W10} for stochastic differential equations (SDEs) with multiplicative noise, and in \cite{ES} for stochastic differential delay equations (SDDEs) with additive noise.
 The aim of this paper is to extend these existed results to the functional solution of SDDEs with multiplicative noise.   Due to the the double difficulty caused by  delay and  non-constant diffusion coefficient, both couplings constructed in \cite{W10} and \cite{ES} are no longer valid. Under a reasonable assumption (see {\bf (A)} below), we will construct a successful coupling which leads to an explicit log-Harnack inequality of the functional solution (see Theorem \ref{T1.1} below). This weaker version of Harnack inequality was introduced in  \cite{RW, W10b} for elliptic diffusion processes, and it is powerful enough to imply some regularity properties of the semigroup such as the strong Feller property and heat kernel estimates w.r.t. quasi-invariant  probability measures (see Corollary \ref{C1.2} below). The dimension-free Harnack inequality in the sense of \cite{W97} is also derived (see Theorem \ref{T4.1} below).

\

Let $r_0>0$ be fixed, and let $\C=C([-r_0,0];\R^d)$ be equipped with the uniform norm $\|\cdot\|_\infty$. Let $\B_b(\C)$ be the set of all bounded measurable functions on $\C$. Let $B(t)$ be a $d$-dimensional Brownian motion on a complete filtered probability space $(\OO,\{\F_t\}_{t\ge 0}, \P)$, and let
\beg{equation*}\beg{split}
&  \si:   [0,\infty)\times\R^d\times \OO\to \R^d\otimes \R^d,\\
&Z:  [0,\infty)\times\R^d\times \OO\to \R^d,\\
&b:  [0,\infty)\times\C\times \OO\to \R^d\end{split}\end{equation*} be progressively measurable and continuous w.r.t. the second variable.  Consider the following delay SDE on $\R^d$:

\beq\label{1.1} \d X(t)= \big\{Z(t,X(t))+b(t,X_t)\big\}\d t +\si(t, X(t))\d B(t),\ \ X_0\in \C,\end{equation} where for each $t\ge 0$, $X_t\in \C$ is fixed as $X_t(u)=X(t+u), u\in [-r_0, 0].$ Let $\|\cdot\|$ and $\|\cdot\|_{HS}$ be the operator norm and the Hilbert-Schmidt norm for $d\times d$-matrices respectively.

To ensure the existence, uniqueness, non-explosion, and further regular properties of the solution, we make use of the following assumption:

\paragraph{(A)} \ \emph{  $\si$ is invertible, and there exist constants $K_1, K_2\ge 0, K_3>0$ and $K_4\in \R$ such that \beg{enumerate}
\item[$(A1)$] $\big|\si(t, \eta(0))^{-1} \{b(t,\xi)-b(t,\eta)\}\big|\le K_1 \|\xi-\eta\|_\infty,\ t\ge 0, \xi,\eta\in \C;$\item[$(A2)$] $\big|(\si(t,x)-\si(t,y))\big|\le K_2 (1\land |x-y|),\  t\ge 0, x,y\in \R^d;$
    \item[$(A3)$] $\big|\si(t,x)^{-1} \big|\le K_3,\  t\ge 0, x\in\R^d;$ \item[$(A4)$] $\big\|\si(t,x)-\si(t,y)\|_{HS}^2 + 2\<x-y, Z(t,x)-Z(t,y)\>\le K_4 |x-y|^2,\ t\ge 0, x,y\in \R^d$\end{enumerate} hold almost surely. }

        \

We remark that in \cite{ES}  $\si$ is assumed to be the unit matrix and $(A4)$ holds for non-positive $K_4$, so that {\bf (A)}  holds for $K_2=0, K_3=1$ and $K_1$ being the Lipschitz constant of $b$.  Moreover, it is easy to see that
 {\bf (A)} is satisfied provided $\si$ is uniformly invertible, and $\si, b, Z$ are Lipschitz continuous w.r.t. the second variable uniformly in the first and third variables. Therefore, our framework is much more general.

On the other hand, if {\bf (A)} holds then for any $\F_0$-measurable $X_0$,
the equation (\ref{1.1}) has a unique strong solution and the solution is non-explosive.
To indicate the dependence of the solution on the initial data, for any $\xi\in\C$ we shall use $X^\xi(t)$ and $X_t^\xi$  respectively to denote the solution and the functional solution  to the equation with $X_0=\xi.$ We shall investigate the Harnack inequality and applications for the family of Markov operators $(P_t)_{t\ge 0}$ on $\B_b(\C)$ given by
$$P_t f(\xi)= \E f(X_t^\xi),\ \ t\ge 0, f\in \B_b(\C).$$
We note that due to the delay, the solution $X(t)$ is not Markovian. But when $Z, b$ and $\si$ are deterministic,   the   functional solution $X_t$ is a strong Markov process.

\beg{thm}\label{T1.1} Assume {\bf (A)}. Then the log-Harnack inequality
$$ P_T\log f(\eta)\le  \log P_T f(\xi) + H_T(\xi,\eta),\ \ T>r_0, \xi,\eta\in\C$$
holds for   $f\in \B_b(\C)$ with $f\ge 1$ and
 $$H_T(\xi,\eta):= \inf_{s\in (0, T-r_0]} \bigg\{\ff{2K_3^2K_4|\xi(0)-\eta(0)|^2}{1-\e^{-K_4 s}}
 +K_1^2\Big\{\ff {r_0}2 +s\big(1+K_2^2K_3^2\big)\Big\} \e^{K_2^2(K_1^2s+8)s} \|\xi-\eta\|_\infty^2\bigg\}.$$
 Consequently, for any $T>r_0$, $P_T$ is strong Feller, i.e. $P_T\B_b(\C)\subset C_b(\C),$ the set of bounded continuous functions on $\C$.\end{thm}

It is easy to see that the log-Harnack inequality only holds for $T>r_0$. Indeed, if the inequality holds for some $T\in (0,r_0]$ then by taking
$f(\xi)=(1+|\xi(T-r_0)|\land n)^n$ and letting $n\to\infty$, the inequality implies that
$$\log(1+|\eta(T-r_0)|) \le \log (1+|\xi(T-r_0)|)$$ holds for all $\xi,\eta\in\C$, which is however impossible.

\

Next, we present some consequences of the above log-Harnack inequality for heat kernels of $P_T$ w.r.t. a quasi-invariant probability measure $\mu$.

\beg{defn} \emph{Let $(E,\F)$ be a measurable space with $\B_b(E)$ the set of all bounded measurable functions,   let $\mu$ be a probability measure on $E$, and let $P$ be a bounded linear operator on $\B_b(E)$. \beg{enumerate} \item[$(i)$]    $\mu$   is called   quasi-invariant   of $P$, if
 $ \mu P$ is absolutely continuous w.r.t. $\mu$, where $(\mu P)(A):= \mu(P1_A),\ A\in\F$. If   $\mu P=\mu$ then $\mu$ is called an invariant probability measure of $P$. \item[$(ii)$]  A   measurable function $p$ on $E^2$ is called the   kernel of  $P$   w.r.t. $\mu$, if
  $$Pf= \int_E p(\cdot, y)f(\eta)\mu(\d y),\ \ \ f\in \B_b(E).$$  \end{enumerate}}\end{defn}

\beg{cor} \label{C1.2} Assume {\bf (A)}. Let $t>r_0$ and $\mu$ be a quasi-invariant probability measure of $P_t$.  Then:
\beg{enumerate} \item[$(1)$]   $P_t$ has a  kernel $p_t$ w.r.t. $\mu$. \item[$(2)$]  The kernel $p_t$ satisfies the entropy inequality
$$\int_\C p_t(\xi,\cdot)\log \ff{p_t(\xi,\cdot)}{p_t(\eta,\cdot)}\,\d\mu
  \le H_t(\xi,\eta),\ \ \xi,\eta\in\C,$$ where we set $r\log \ff r s =0$ if $r=0$ and $r\log \ff r s =\infty$ if $r>0$ and $s=0$.
\item[$(3)$]  The kernel $p_t$ satisfies
$$\int_\C p_t(\xi,\cdot) p_t(\eta,\cdot)\,\d\mu \ge \exp[- H_t(\xi,\eta)],\ \ \xi,\eta\in\C.$$ \item[$(4)$] $P_t$   has at most one invariant probability measure, and if it has,   the kernel of $P_t$ w.r.t. the invariant probability measure is strictly positive.
\end{enumerate}\end{cor}
Note that if $P_t$ is symmetric w.r.t. $\mu$, then $\int_\C p_t(\xi,\cdot)p_t(\eta,\cdot)\d\mu=p_{2t}(\xi,\eta)$ so that (3) provides a Gaussian type lower bound for the heat kernel. Moreover, if $\mu$ is an invariant probability measure of $P_t$, then (2) gives an entropy-cost inequality as in \cite[Corollary 1.2(3)]{RW}. More precisely, letting $P_t^*$ be the adjoint operator of $P_t$ in $L^2(\mu)$,
for any $f\ge 0$ with $\mu(f)=1$, one has
$$\mu((P_t^*f)\log P_t^*f) \le \inf_{\pi\in \mathcal C(\mu,f\mu)} \int_{\C\times \C} H_t(\xi,\eta)\pi(\d\xi,\d\eta),\ \ t>0,$$ where $\mathcal C(\mu, f\mu)$ is the set of all couplings for $\mu$ and $f\mu$. The right hand side of the above inequality is called the transportation-cost between $\mu$ and $f\mu$ with cost function $H_t$.
Finally, we note that the uniqueness of the invariant probability measure has been investigated in \cite{MH} for SDDEs in terms of the asymptotic coupling property.

\

To conclude this section, let us present an existence result of the quasi-invariant measure (see \cite{EGS09} for existence of the invariant probability measure).

\beg{prp}\label{P1.3} Assume {\bf (A)} and let $\si$ and $Z$   be deterministic and time-independent. If  $x\mapsto -\|\si(x)\|_{HS}^2 -\<x, Z(x)\>$ is a compact function, i.e. $\{x\in\R^d: -\|\si(x)\|_{HS}^2- \<x,Z(x)\>\le r\}$ is a compact set for any constant $r>0$, then $\{P_t\}_{t\ge 0}$ has a  quasi-invariant probability measure, i.e. the measure is quasi-invariant for all $P_t, t\ge 0$. \end{prp}

The proof of Theorem \ref{T1.1} is presented in Section 2 while those of Corollary \ref{C1.2} and Proposition \ref{P1.3} are addressed in Section 3. Finally, in Section 4 we investigate the dimension-free Harnack inequality in the sense of \cite{W97}.

\section{Proof of Theorem 1.1}

According to \cite[Proposition 2.3]{W10}, the claimed log-Harnack inequality implies the strong Feller property of $P_T$, see also Proposition \ref{P}(1) below. So, we only have to prove the desired log-Harnack inequality.
To make the proof easy to follow, let us first explain the main idea of the argument.

\

Let $T>r_0$ and $t_0\in (0, T-r_0]$ be fixed.
Let $X(s)$ solve (\ref{1.1}) with $X_0=\xi$. For $\gg\in C^1([0,t_0])$ such that $\gg(r)>0$ for $r\in [0,t_0)$ and $\gg(t_0)=0$, let $Y(t)$ solve the equation
\beq\label{2.1} \beg{split}  \d Y(t)= &\Big\{Z(t, Y(t))+b(t,X_t) + \ff{1_{\{t<t_0\}}}{\gg(t)}\si(t,Y(t)) \si(t, X(t))^{-1} (X(t)-Y(t))\Big\}\d t\\
&\qquad \qquad+\si(t, Y(t))\d B(t),\ \ Y_0=\eta.\end{split}\end{equation}
The key point of our coupling is that  $X(t)$ and $Y(t)$ will move together from time $t_0$ on, so that $X_T=Y_T$. To this end, we add the  drift term $$\ff{1_{\{t<t_0\}}}{\gg(t)} \si(t,Y(t))\si(t, X(t))^{-1} (X(t)-Y(t)) \d t$$ to force $Y(t)$ to meet $X(t)$ at time $t_0$. In order to dominate the non-trivial martingale part of $X(t)-Y(t)$,   the force has to be infinitely strong near by $t_0$, for this we need $\gg(t_0)=0$. More precisely, as in \cite{W10} we shall take
$$\gg(t) = \ff{2-\theta}{K_4} \big(1-\e^{(t-t_0)K_4}\big),\ \ t\in [0,t_0)$$ for a parameter $\theta\in (0,2)$. In this case, we have
\beq\label{GG} 2 + \gg'(t)-K_4\gg(t)= \theta,\ \ t\in [0, t_0].\end{equation}

Moreover, to ensure   these two process moving together after the coupling time (i.e. the first meeting time), they should solve the same equation  from that time on. This is the reason why we have to take the delay term in (\ref{2.1}) by using $X_t$ rather than $Y_t$.  Since the additional drift is singular at time $t_0$, it is only clear that $Y(t)$ is well solved before time $t_0$. To solve $Y(t)$ for all $t\in [0,T]$, we need to reformulate the equation by using a new Brownian motion determined by the Girsanov transform induced by the coupling.

\

Let
$$\phi_t= \si(t,Y(t))^{-1} \{b(t, Y_t)-b(t,X_t)\} -\ff {1_{\{t<t_0\}}} {\gg(t)} \si(t, X(t))^{-1} (X(t)-Y(t)),\ \ t\ge 0.$$ From {\bf (A)} it is easy to see that

$$R_t:= \exp\bigg[\int_0^t\<\phi_s, \d B(s)\> -\ff 1 2\int_0^t |\phi_s|^2\d s\bigg]$$ is a martingale for $t\in [0,t_0)$. We shall further prove that

\ \newline {\bf \ (i)} $ \{R_t\}_{t\ge 0}$ is a well-defined  martingale.

\ \newline Whence {\bf (i)} is confirmed, by the Girsanov theorem, under probability $\d \Q_{T}:= \R_{T}\d\P$ the process
$$\tt B_t:= B_t -\int_0^t \<\phi_s, \d B(s)\>,\ \ t\in [0,T]$$ is a $d$-dimensional Brownian motion  and (\ref{2.1}) reduces to
 \beq\label{*2}   \d Y(t)= \big\{Z(t, Y(t))+b(t,Y_t)\big\}\d t+\si(t, Y(t))\d \tt B(t),\ \ t\in [0,T], Y_0=\eta.\end{equation}
Therefore, (\ref{2.1}) has a unique solution $\{Y(t)\}_{t\in [0,T]}$ under the probability $\Q_T$, and
\beq\label{C1}P_Tf(\eta)= \E_{\Q_T} f(Y_T)= \E[R_T f(Y_T)]. \end{equation}
Next, we shall prove that

\ \newline
 {\bf (ii)}  The coupling time $\tau:=\inf\{t\ge 0: X(t)=Y(t)\}\le t_0,\ \Q_T$-a.s.

 \

From (\ref{1.1}) and (\ref{2.1}) we see that for $t\ge \tau$, the two processes  $X(t)$ and $Y(t)$   solve the same equation,
because  the additional drift term   disappears as soon as $X(t)=Y(t)$. By  the uniqueness of the solution to (\ref{1.1})
we have $X(t)=Y(t)$ for $t\ge \tau$. Combining this with {\bf (ii)} and noting that $t_0\le T-r_0$,  we conclude that   $X_T=Y_T,\ \Q_T$-a.s.  So, by the Young inequality and (\ref{C1}), we arrive at
$$ P_T \log f(\eta)= \E[R_T \log f(Y_T)]= \E[R_T\log f(X_T)] \le \log P_T f(\xi)+ \E R_T\log R_T.$$ Therefore, to complete the proof it remains to show that

\ \newline {\bf (iii)} $\E \big[R_T\log R_T\big]\le \ff{2K_3^2K_4|\xi(0)-\eta(0)|^2}{1-\e^{-K_4t_0}}
 + K_1^2\Big\{\ff {r_0}2 +t_0(1+K_2^2K_3^2)\Big\} \e^{K_2^2(K_1^2t_0+8)t_0} \|\xi-\eta\|_\infty^2.$

\

In the remainder of the section, we will prove the above claimed {\bf (i)}-{\bf (iii)} respectively.

\subsection{Proofs of (i)}
 The key result of this subsection is the following.

\beg{prp}\label{P2.1} Assume {\bf (A)}. Then for any $t\in [0, t_0)$,
$$\E \big[R_t\log R_t \big]\le \ff{2K_3^2K_4|\xi(0)-\eta(0)|^2}{\theta(2-\theta)(1-\e^{-K_4t_0})} +\ff{tK_1^2(1 +K_2^2K_3^2) \e^{K_2^2(K_1^2t+8)t}}{\theta^2}\|\xi-\eta\|_\infty^2.$$
 \end{prp}
\beg{proof}
Let $t\in (0,t_0)$ be fixed. Then $\{\tt B(s)\}_{s\le t}$ is a $d$-dimensional Brownian motion under the probability $\d \Q_t:=R_t\d \P,$ so that $(A1)$ and $(A3)$ imply
\beq\label{**} \beg{split} \E [R_t\log R_t]& =\E_{\Q_t} \log R_t = \E_{\Q_t} \bigg\{\int_0^t \<\phi_s, \d\tt B(s)\>
+\ff 1 2 \int_0^t |\phi_s|^2\d s\bigg\}\\
&=\ff 1 2 \int_0^t \E_{\Q_t} |\phi_s|^2\d s\\
&\le K_1^2\int_0^t \E_{\Q_t} \|Y_s-X_s\|_\infty^2\d s + K_3^2\int_0^t\ff 1 {\gg(s)^2} \E_{\Q_t} |X(s)-Y(s)|^2\d s\\
&=:I_1+I_2.\end{split}\end{equation} To estimate $I_1$ and $I_2$, let us reformulate equation (\ref{1.1}) using the new Brownian motion $\tt B(s)$:
$$\d X(s) =\big\{Z(s, X(s))+ b(s, X_s)+ \si(s, X(s))\phi_s\big\}\d s +\si(s, X(s))\d\tt B(s),\ \ s\le t.$$ Since \beg{equation*} \beg{split} &\si(s,X(s))\phi_s +b(s, X_s)- b(s, Y_s) \\
&=\big\{\si(s, X(s))-\si(s, Y(s))\big\}\si(s, Y(s))^{-1} (b(s, Y_s)-b(s, X_s))
-\ff{X(s)-Y(s)}{\gg(s)}, \end{split}\end{equation*} the equation reduces to
\beq\label{*3}\beg{split} \d X(s)= \bigg\{&\{\si(s, X(s))-\si(s, Y(s))\}\si(s, Y(s))^{-1} (b(s, Y_s)-b(s, X_s))\\
& +Z(s, X(s)) + b(s, Y_s)-\ff{X(s)-Y(s)}{\gg(s)}\bigg\}\d s +\si(s, X(s))\d \tt B(s),\ \ s\le t.\end{split}\end{equation}Combining this with (\ref{*2}) and using the It\^o formula, we obtain from $(A1), (A2)$ and $(A4)$ that
\beq\label{*4}\beg{split}&\d |X(s)-Y(s)|^2  \le 2 \big\<X(s)-Y(s), (\si(s, X(s))-\si(s, Y(s)))\d\tt B(s)\big\> \\
&+\bigg\{2K_1K_2 \|X_s-Y_s\|_\infty|X(s)-Y(s)| +\Big(K_4-\ff 2 {\gg(s)} \Big)|X(s)-Y(s)|^2 \bigg\}\d s,\ s\le t.\end{split}\end{equation}
Since it is easy to see that
$$K_4\le \ff 2 {\gg(0)}\le \ff 2 {\gg (s)},$$  it follows that
\beq\label{*5} \beg{split}\d |X(s)-Y(s)|\le &\Big\<\ff{X(s)-Y(s)}{|X(s)-Y(s)|}, (\si(s, X(s))-\si(s, Y(s)))\d\tt B(s)\Big\> \\
&+K_1K_2 \|X_s-Y_s\|_\infty \d s,\ \ s\le t.\end{split}\end{equation} Let
$$M(s)= \int_0^s \Big\<\ff{X(r)-Y(r)}{|X(r)-Y(r)|}, (\si(r, X(r))-\si(r, Y(r)))\d\tt B(r)\Big\>,\ \ s\le t,$$ which is a martingale under $\Q_t$. By $(A2)$ and the Doob inequality we have
$$\E_{\Q_t} \sup_{r\in [0,s]} M(r)^2\le 4 K_2^2 \int_0^s \E_{\Q_t}\|X_r-Y_r\|_\infty^2\d r,\ \ s\le t.$$ Combining this with (\ref{*5}) we obtain
$$\E_{\Q_t} \|X_s-Y_s\|_\infty^2 \le \|\xi-\eta\|_\infty^2 + 2K_2^2(K_1^2s+4) \int_0^s \|X_r-Y_r\|_\infty^2 \d r,\ \ s\le t.$$ By the Gronwall lemma, this  implies that
\beq\label{*6} \E_{\Q_t} \|X_s-Y_s\|_\infty^2\le \|\xi-\eta\|_\infty^2 \e^{K_2^2(K_1^2s+8)s},\ \ s\le t.\end{equation} On the other hand, let
$$\d \tt M(s)= \ff 2 {\gg(s)} \big\<X(s)-Y(s), (\si(s,X(s))-\si(s,Y(s)))\d\tt B(s)\big\>,\ \ s\le t,$$ which is a martingale under $\Q_t$. It follows from (\ref{*4}) and (\ref{GG}) that
\beg{equation}\label{*7}\beg{split} &\d\ff{|X(s)-Y(s)|^2}{\gg(s)}-\d\tt M(s)\\
&\le \bigg\{\ff{K_1K_2}{\gg(s)} \|X_s-Y_s\|_\infty|X(s)-Y(s)| +\ff{K_4 \gg(s)-2-\gg'(s)}{\gg(s)^2}|X(s)-Y(s)|^2\bigg\}\d s\\
&\le \bigg\{\ff {K_1K_2}{\gg(s)}  \|X_s-Y_s\|_\infty|X(s)-Y(s)|  -\ff{\theta |X(s)-Y(s)|^2}{\gg(s)^2}\bigg\}\d s, \ \ s\le t.\end{split}\end{equation}
Combining this with   (\ref{*6}), we obtain
\beg{equation*}\beg{split}& h(t):=  \int_0^t \ff{\E_{\Q_t}|X(s)-Y(s)|^2}{\gg(s)^2} \d s\\
&\le \ff{|\xi(0)-\eta(0)|^2}{\theta\gg(0)} +\ff{K_1K_2}{\theta}h(t)^{1/2} \bigg(\int_0^t \E_{Q_t} \|X_s-Y_s\|_\infty^2\d s\bigg)^{1/2} \\
&\le  \ff{|\xi(0)-\eta(0)|^2}{\theta\gg(0)} + \ff {h(t)} 2 + \ff{K_1^2K_2^2} {2 \theta^2} \e^{K_2^2(K_1^2t+8)t}\|\xi-\eta\|_\infty^2.\end{split}\end{equation*} Therefore,
\beq\label{INT}\int_0^t \ff{\E_{\Q_t}|X(s)-Y(s)|^2}{\gg(s)^2} \d s\le  \ff{2|\xi(0)-\eta(0)|^2}{\theta\gg(0)}  +  \ff{K_1^2K_2^2t  \e^{K_2^2(K_1^2t+8)t}}{\theta^2}\|\xi-\eta\|_\infty^2.\end{equation} Substituting this and (\ref{*6}) into (\ref{**}), we complete the proof.\end{proof}

\beg{proof}[Proof of {\bf (i)}] According to Proposition \ref{P2.1}, $\{R_t\}_{t\in [0, t_0)}$ is a uniformly integrable continuous martingale. So, by the martingale convergence theorem, \beq\label{2.**}R_{t_0}=\lim_{t\uparrow t_0} R_t\end{equation}
exists and $\{R_t\}_{t\in [0, t_0]}$ is again a uniformly integrable martingale. In particular, $\tt B(t)$ is a $d$-dimensional Brownian motion under $\d\Q_{t_0}:= R_{t_0}\d\P$ such that $Y(t)$ can be solved from (\ref{*2}) for $t\in [0,t_0]$. To solve the equation for $t> t_0$, let us use the filtered probability space
$(\OO, \{\F_t\}_{t\ge t_0}, \Q_{t_0})$. Since $\{B(t)-B(t_0)\}_{t\ge t_0}$ is independent of $R_{t_0}$, it is easy to see that $\{B(t)\}_{t\ge t_0}$ is  a $d$-dimensional Brownian motion on this probability space. Moreover, due to (\ref{*6}),  $X_{t_0}$ and $Y_{t_0}$ are $\F_{t_0}$-measurable random variables on $\C$ with

\beq\label{*6'} \E_{\Q_{t_0}} \|X_{t_0}-Y_{t_0}\|_\infty^2\le \|\xi-\eta\|_\infty^2
\e^{K_2^2(K_1^2t_0+8)t_0}.\end{equation} Therefore, by $(A2)$ and $(A4)$, starting from $Y_{t_0}$ at time $t_0$ the equation (\ref{2.1}) has a unique solution $\{Y(t)\}_{t\ge t_0}$, and by the It\^o formula and $(A4)$,

$$\d|X(t)-Y(t)|^2\le 2\<X(t)-Y(t), (\si(t,X(t))-\si(t,Y(t)))\d B(t)\> +K_4 |X(t)-Y(t)|^2\d t,\ \ t\ge t_0.$$
Combining this with $(A2)$ and noting that $\{B(t)\}_{t\ge t_0}$ is a $\Q_{t_0}$-Brownian motion, we obtain

$$\E_{\Q_{t_0}}\Big(\sup_{s\in [t_0, t]}|X(s)-Y(s)|^2\Big|\F_{t_0}\Big)\le \e^{K(t-t_0)}\|X_{t_0}-Y_{t_0}\|_\infty^2,\ \ t\ge t_0$$ for some constant $K>0$. Therefore, it follows from (\ref{*6'})  that

$$\E_{\Q_{t_0}}\sup_{s\in [t_0,t]} \|X_s-Y_s\|_\infty^2<\infty,\ \ t\ge t_0.$$ Since

$$|\phi_t|=\big|\si(t,Y(t))^{-1}(b(t,Y_t)-b(t,X_t))\big|\le K_1 \|X_t-Y_t\|_\infty,\ \ t\ge t_0,$$ this implies that
$$\ff{R_t}{R_{t_0}}=\exp\bigg[\int_{T-r_r}^t\<\phi_s,\d B(s)\>- \ff 1 2 \int_{t_0}^t |\phi_s|^2\d s\bigg],\ \ t\ge t_0$$ is a $\Q_{t_0}$-martingale, and thus, for $t>s\ge t_0$ and $A\in \F_s$,

$$\E(R_t1_A)= \E_{\Q_{t_0}} \Big\{1_A\ff{R_t}{R_{t_0}}\Big\}= \E_{\Q_{t_0}}\Big\{1_A\ff{R_s}{R_{t_0}}\Big\}=
\E(R_s1_A).$$ This means that $\{R_t\}_{t\ge t_0}$ is a $\P$-martingale, and thus,    $\{R_t\}_{t\ge 0}$ is a well-defined $\P$-martingale as claimed since $\{R_t\}_{t\in [0,t_0]}$ is already a martingale.
\end{proof}

\subsection{Proof of (ii)} Since $\{R_t\}_{t\in [0,T]}$ is a martingale,  for any $t\in [0,t_0)$ the inequality (\ref{INT}) holds for  $\Q_T$ in place of  $\Q_t$.   Therefore,
\beq\label{N}\E_{\Q_{T}} \int_0^{t_0}\ff{|X(t)-Y(t)|^2}{\gg(t)^2}\d t<\infty.\end{equation} This implies that $\tau\le t_0,\ \Q_T$-a.s.
Indeed, since  $t\mapsto X(t)$ and $t\mapsto Y(t)$ are continuous  $\Q_T$-a.s., there exists $\OO_0\subset \OO$ with $\Q_T(\OO_0)=1$ such that for any $\oo\in\OO_0$, $X(t)(\oo)$ and $Y(t)(\oo)$ are continuous in $t$.  If $\oo\in \OO_0$ such that $\tau(\oo)>t_0$, then

$$\inf_{t\in [0,t_0]}|X(t)-Y(t)|(\oo)>0,$$ so that

$$\int_0^{t_0}\ff{|X(t)-Y(t)|^2(\oo)}{\gg(t)^2}\d t\ge \inf_{t\in [0,t_0]}|X(t)-Y(t)|(\oo)\int_0^{t_0}\ff {\d t} {\gg(t)^2} =\infty.$$ This means that

$$\Q_T(\tau>t_0) \le \Q_T\bigg(\int_0^{t_0}\ff{|X(t)-Y(t)|^2}{\gg(t)^2}\d t=\infty\bigg)$$ which equals to zero according to (\ref{N}).

\subsection{Proof of (iii)} Since $\{R_t\}_{t\ge 0}$ is a martingale, by the Girsanov theorem   $\{\tt B(t)\}_{t\in [0,T]}$ is Brownian motion under $\Q_T$.  Then

\beq\label{AA} \beg{split} \E [R_T\log R_T] &=\ff 1 2 \E_{\Q_T} \int_0^T |\phi_t|^2\d t \\
& = \ff 1 2 \E_{\Q_{t_0}} \int_0^{t_0}|\phi_t|^2\d t+ \ff 1 2 \E_{\Q_T} \int_{t_0}^T |\phi_t|^2\d t\\
 &= \E [R_{t_0}\log R_{t_0}] + \ff 1 2 \E_{\Q_T} \int_{t_0}^T |\phi_t|^2\d t.\end{split}\end{equation} Since $\tau\le t_0$ and $X(t)=Y(t)$ for $t\ge \tau$, we have $X(t)=Y(t)$ for  $t\ge t_0$. So, it follows form $(A1)$ that

 $$\int_{t_0}^T |\phi_t|^2\d t \le K_1^2 \int_{t_0}^T \|X_t-Y_t\|_\infty^2 \d t\le K_1^2 r_0 \|Y_{t_0}-X_{t_0}\|_\infty^2.$$ Combining this with (\ref{*6}),  which also  holds  for $t=s=t_0$ by (\ref{2.**}) and the Fatou lemma, we arrive at

 $$\ff 1 2 \E_{\Q_T} \int_{t_0}^T |\phi_t|^2\d t \le \ff {K_1^2 r_0} 2\e^{K_2^2(K_1^2t_0+8)t_0} \|\xi-\eta\|_\infty^2.$$ Substituting this into (\ref{AA}) and noting that (\ref{2.**}) and  Proposition \ref{P2.1} with $\theta=1$ imply
$$ \E [R_{t_0}\log R_{t_0}]\le \ff{2K_3^2K_4|\xi(0)-\eta(0)|^2}{\theta(2-\theta)(1-\e^{-K_4t_0})} +\ff{t_0K_1^2(1 +K_2^2K_3^2) \e^{K_2^2(K_1^2t_0+8)t_0}}{\theta^2}\|\xi-\eta\|_\infty^2,$$
 we prove {\bf (iii)}.

 \section{Proofs of Corollary \ref{C1.2} and Proposition \ref{P1.3}}

 \beg{proof}[Proof of Corollary \ref{C1.2}]
According to Theorem \ref{T1.1}, we have

$$\e^{P_Tf}(\xi)\le  \{P_T \e^f (\eta)\}\e^{\Psi_T(\|\xi-\eta\|_\infty)},\ \ \xi,\eta\in \C, T>r_0, f\in \B_b(\C), f\ge 0$$ holds for
some continuous function $\Psi_T$ with $\lim_{r \to 0} \Psi_T(r)=0.$
So, the desired assertions follow immediately from the following more general result for $\Phi(r)=\e^r$.
 \end{proof}

 \beg{prp}\label{P} Let $(E,\F)$ be the Borel measurable space of a topology space $E$, $P$ a Markov operator on $\B_b(E)$, and $\mu$ a quasi-invariant probability measure of $P$. Let $\Phi\in C^1([0,\infty))$ be an increasing function with $\Phi'(1)>0$ and
 $\Phi(\infty):=\lim_{r\to\infty}\Phi(r)=\infty$, such that

 \beq\label{LP} \Phi(P f(x))\le  \{P\Phi(f)(y)\}\e^{ \Psi(x,y)},\ \ x,y\in E, f\in \B_b(E), f\ge 0\end{equation} holds for some measurable non-negative function $\Psi$ on $E^2$.
 \beg{enumerate}\item[$(1)$] If $\lim_{y\to x} \{\Psi (x,y)+\Psi(y,x)\}=0$ holds for all $x\in E$, then $P$ is strong Feller. \item[$(2)$] $P$ has a kernel $p$ w.r.t. $\mu$, so that any invariant probability measure of $P$ is absolutely continuous w.r.t. $\mu$.
 \item[$(3)$]   $P$ has at most one invariant probability measure and if it has, the kernel of $P$ w.r.t. the invariant probability measure is strictly positive.
 \item[$(4)$] The kernel $p$ of $P$ w.r.t. $\mu$ satisfies
 $$\int_E p(x,\cdot)\Phi^{-1} \Big(\ff{p(x,\cdot)}{p(y,\cdot)}\Big)\d\mu\le \Phi^{-1}(\e^{\Psi(x,y)}),\ \ x,y\in E,$$
 where $\Phi^{-1}(\infty):=\infty$ by convention.  \item[$(5)$] If $r\Phi^{-1}(r)$ is convex for $r\ge 0$, then the kernel $p$ of $P$ w.r.t. $\mu$ satisfies
  $$\int_E p(x,\cdot)p(y,\cdot)\d\mu\ge \e^{-\Psi(x,y)},\ \ x,y\in E.$$\end{enumerate} \end{prp}

\beg{proof} (1) Let $f\in \B_b(E)$ be positive. Applying (\ref{LP}) to $1+\vv f$ in place of $f$ for $\vv>0$, we have

$$\Phi(1+\vv Pf(x) ) \le \{P\Phi(1+\vv f)(y)\}\e^{\Psi(x,y)},\ \ x,y\in E, \vv>0.$$ By the Taylor expansion this implies
\beq\label {D}\Phi(1)+ \vv \Phi'(1) Pf(x)+ \circ(\vv)\le \{\Phi(1)+ \vv \Phi'(1) Pf(y) +\circ(\vv)\}\e^{\Psi(x,y)}\end{equation} for small $\vv>0$.
Letting $y\to x$ we obtain
$$ \vv Pf(x)\le \vv\liminf_{y\to x} Pf(y)+\circ(\vv).$$ Thus, $Pf(x)\le \lim_{y\to x}Pf(y)$ holds for all $x\in E$. On the other hand, letting $x\to y$ in (\ref{D}) gives $Pf(y)\ge \limsup_{x\to y}Pf(x)$ for any $y\in E$. Therefore, $Pf$ is continuous.

(2) To prove the existence of kernel, it suffices to prove that for any $A\in \F$ with $\mu(A)=0$ we have $P1_A\equiv 0$. Applying (\ref{LP}) to $f= 1 +n 1_A$, we obtain
\beq\label{Z1} \Phi(1+ n P1_A(x))\int_E\e^{-\Psi(x,y)}\mu(\d y)\le \int_E\Phi(1+n1_A)(y)(\mu P)(\d y),\ \ n\ge 1.\end{equation} Since $\mu(A)=0$ and $\mu$ is quasi-invariant for $P$, we have $1_A=0, \mu P$-a.s. So, it follows from (\ref{Z1}) that
$$\Phi(1+n P1_A(x))\le \ff{\Phi(1)}{\int_E\e^{-\Psi(x,y)}\mu(\d y)}<\infty,\ \ x\in E, n\ge 1.$$ Since $\Phi(1+n)\to\infty$ as $n\to\infty$, this implies that $P1_A(x)=0$ for all $x\in E$.

 Now, for any invariant probability measure $\mu_0$ of $P$, if $\mu(A)=0$ then $P1_A\equiv 0$ implies that $\mu_0(A)=\mu_0(P1_A)=0$. Therefore, $\mu_0$ is absolutely continuous w.r.t. $\mu$.

(3)  We first prove that the kernel of $P$ w.r.t. an invariant probability measure $\mu_0$ is strictly positive. To this end, it suffices to show that for any $x\in E$ and $A\in \F$, $P1_A(x)=0$ implies that $\mu_0(A)=0$. Since $P1_A(x)=0$, applying (\ref{LP}) to $f= 1+n P1_A$ we obtain
$$\Phi(1+n P 1_A(y)) \le \{P\Phi(1+n1_A)(x)\}\e^{\Psi(y,x)}=  \Phi(1)\e^{\Psi(y,x)},\ \ y\in E, n\ge 1.$$ Letting $n\to\infty$ we conclude that $P1_A\equiv 0$ and hence, $\mu_0(A)=\mu_0(P1_A)=0.$

Next, let $\mu_1$ be another invariant probability measure of $P$, by (2) we have $\d\mu_1 =f\d\mu_0$ for some probability density function $f$. We aim to prove that $f=1, \mu_0$-a.e. Let $p(x,y)>0$ be the kernel of $P$ w.r.t. $\mu_0$, and let
$P^*(x,\d y)= p(y,x)\mu_0(\d y).$ Then

$$P^* g=\int_E g(y) P^*(\cdot,\d  y),\ \ g\in \B_b(E)$$ is the adjoint operator of $P$ w.r.t. $\mu_0$. Since $\mu_0$ is $P$-invariant, we have

$$\int_E g P^*1\,\d\mu_0= \int_E Pg \,\d\mu_0= \int_E g\,\d\mu_0,\ \ g\in \B_b(E).$$ This implies that $P^*1=1, \mu_0$-a.e. Thus, for $\mu_0$-a.e. $x\in E$ the measure $P^*(x,\cdot)$ is a probability measure. On the other hand, since $\mu_1$ is $P$-invariant, we have

$$\int_E(P^*f)g\, \d\mu_0 =\int_E f Pg\, \d\mu_0 = \int_E Pg\, \d\mu_1= \int_E g\,\d\mu_1= \int_E fg\,\d\mu_0,\ \ g\in \B_b(E).$$
This implies that $P^*f=f, \mu$-a.e. Therefore, for any $r>0$ we have

$$\int_E P^*\ff 1 {f+1}\,\d\mu_0= \int_E \ff 1 {f+1} \,\d\mu_0 =\int_E\ff 1 {P^*f+1}\,\d\mu_0.$$ When $P^*(x,\cdot)$ is a probability measure, by the Jensen inequality   one has $P^*\ff 1 {1+f}(x)\ge \ff 1 {P^*f+1}(x)$ and the equation holds if and only if  $f$ is constant $P^*(x,\cdot)$-a.s. Hence,   $f$ is constant $P^*(x,\cdot)$-a.s. for $\mu_0$-a.e. $x$. Since $p(x,y)>0$ for any $y\in E$ such that $\mu_0$ is absolutely continuous  w.r.t. $P^*(x,\cdot)$ for any $x\in E$, we conclude that $f$ is constant $\mu_0$-a.s. Therefore,  $f=1\ \mu_0$-a.s. since $f$ is a probability density function.

(4) Applying (\ref{LP}) to
$$f=n\land \Phi^{-1} \Big(\ff{p(x,\cdot)}{p(y,\cdot)}\Big)$$ and letting $n\to\infty$, we obtain the desired inequality.

(5) Let $r\Phi^{-1}(r)$ be convex for $r\ge 0$. By the Jensen inequality we have
$$\int_E p(x,\cdot)\Phi^{-1} (p(x,\cdot))\d\mu\ge \Phi^{-1}(1).$$
So, applying (\ref{LP}) to

$$f=n\land \Phi^{-1}(p(x,\cdot))$$ and letting $n\to\infty$,   we obtain
$$\int_E p(x,\cdot)p(y,\cdot)\d\mu \ge \e^{-\Psi(x,y) } \Phi\bigg(\int_E p(x,\cdot)\Phi^{-1} (p(x,\cdot))\d\mu\bigg)\ge \e^{-\Psi(x,y)}.$$ \end{proof}

 \beg{proof}[Proof of Proposition \ref{P1.3}]  Consider the It\^o SDE without delay:
 \beq\label{F1} \d \tt X(t)= Z(\tt X(t))\d t +\si (\tt X(t))\d B(t).\end{equation} The equation has a unique strong solution which is a strong Markov process. Since $-\|\si(x\|_{HS}^2-\<x,Z(x)\>$ is a compact function, it is standard   that the process has a (indeed unique, due to the ellipticity)
 invariant probability measure $\mu_0$, so that with initial distribution $\mu_0$ the process is stationary. Let $\mu$ be the distribution of the $\C$-valued random variable
 $$\{\tt X(r_0+u)\}_{u\in [-r_0,0]},$$ where $\tt X(0)$ has distribution $\mu_0$. Then the $\C$-valued Markov process $\{X_t\}_{t\ge 0}$ with

 $$X_t(u):= X(t+u):= \tt X(t+r_0+u),\ \ \ u\in [-r_0,0]$$ has an invariant probability measure $\mu$. Let $\bar B(t)= B(t+r_0)-B(r_0)$, we have

 \beq\label{F2} \d X(t)= Z((X(t))\d t +\si(X(t))\d \bar B(t),\ \ t\ge 0.\end{equation} As before, let $X^\xi(t)$ be the solution of this equation with $X_0=\xi$.
 To formulate $P_tf(\xi)$ using $X_t^\xi$, we take

 $$\tt B^\xi(t)=\bar B(t)+\int_0^t \si(X^\xi(s))^{-1} b(s, X_s^\xi) \d s.$$ Then   (\ref{F2}) implies that

 $$\d X^\xi(t)= \{Z(X^\xi(t)) + b(t, X_t^\xi)\}\d t +\si(X^\xi(t))\d\tt B^\xi(t),\ \ t\ge 0.$$  By {\bf (A)} it is easy to see that

 $$R_t^\xi:= \exp\bigg[\int_0^t \<b(s, X_s^\xi), \si(X^\xi(s))\d\bar B(s)\> -\ff 1 2 \int_0^t |\si(X^\xi(s))^{-1} b(s, X_s^\xi)|^2 \d s\bigg]$$ is a martingale, and by the Girsanov theorem for any $T>0$, $\{\tt B^\xi(t)\}_{t\in [0,T]}$ is a Brownian motion under probability $\d\Q_T^\xi:= R_T^\xi\d\P$.  Therefore,

 \beq\label{F3} P_T f(\xi)= \E[R_T^\xi f(X_T^\xi)],\ \ \ T\ge 0, f\in \B_b(\C).\end{equation} Since $\mu$ is an invariant probability measure of $X_t$, for any $\mu$-null set $A$ we have

 $$\int_{\OO\times \C} 1_A(X_T^\xi(\oo))(\P\times\mu)(\d\oo, \d\xi) = \mu(A)=0.$$ Combining this with (\ref{F3}) we obtain

 $$(\mu P_T)(A)= \int_\C P_T 1_A \d\mu = \int_{\OO \times\C}\big\{ 1_A(X_T^\xi(\oo) R_T^\xi(\oo) \big\}(\P\times\mu)(\d\oo, \d\xi)=0.$$ Therefore, $\mu$ is a quasi-invariant probability measure of $P_T$.
 \end{proof}

\section{The Harnack inequality}

 In this section we aim to establish the Harnack inequality with a power $p>1$ in the sense of \cite{W97}:

 \beq\label{4.1} P_T f(\eta)\le \{P_T f^p(\xi)\}^{1/p}\exp[\Phi_p(T,\xi,\eta)],\ \ f\ge 0, T>r_0, \xi,\eta\in\C\end{equation} for some positive function $\Phi_p$ on $(r_0,\infty)\times \C^2.$ As shown in \cite{W10}
 for the case without delay, we will have to assume that $p>(1+K_2K_3)^2.$  In this case, letting
 $$\ll_p= \ff 1 {2(\ss p -1)^2},$$
  the set
 $$\Theta_p:=\bigg\{ \vv\in  (0,1): \ff{(1-\vv)^4}{2(1+\vv)^3K_2^2K_3^2}\ge \ll_p\bigg\}$$ is non-empty.   Let $$W_{ \vv}(\ll)= \max\bigg\{\ff{8(1+\vv)r_0K_1^3K_2\ll\{4(1+\vv)r_0K_1K_2\ll+\vv\}}{\vv^2}, \ff{2(1+\vv)^2\ll}{\vv^2}, \ff{(1+\vv)^3K_1^2K_2^2K_3^2\ll}{8\vv^2(1-\vv)^3}\bigg\},$$ and
 $$s_{\vv}(\ll)= \ff{\ss{K_1^2+2W_{\vv}(\ll)}-K_1}{4W_{\vv}(\ll)K_2},\ \ \vv\in (0,1),\ll>0.$$

\beg{thm}\label{T4.1} Assume {\bf (A)}. For any   $p>(1+K_2K_3)^2$ and $T>r_0$, the Harnack inequality $(\ref{4.1})$ holds for
\beg{equation*}\beg{split} \Phi_p(T,\xi,\eta):=&\ff{\ss p-1}{\ss p}  \inf_{\vv\in\Theta_p}\inf_{s\in (0, s_{\vv}(\ll_p)\land (T-r_0)]} \bigg\{\ff\vv {2(1+\vv)} +\ff{16K_2^2s^2W_{\vv}(\ll_p)}{1-4K_1K_2s}\\ &\quad + \ff{\ll_p(1+\vv)^2K_3^2K_4|\xi(0)-\eta(0)|^2}{2\vv (1-\vv)^2(1+2\vv)(1-\e^{-K_4s})}+ \big(K_1^2r_0\ll_p+2sW_{\vv}(\ll_p)\big)\|\xi-\eta\|_\infty^2\bigg\}.\end{split}\end{equation*}
 Consequently, there exists a decreasing function $C: ((1+K_2K_3)^2,\infty)\to (0,\infty)$ such that $(\ref{4.1})$ holds for
$$\Phi_p(T,\xi,\eta)= C(p)\Big\{ 1+ \ff{|\xi(0)-\eta(0)|^2}{T-r_0}+\|\xi-\eta\|_\infty^2\Big\}.$$ \end{thm}

\beg{proof} (a) We first observe that the second assertion is a consequence of the first. Indeed, for any $q>(1+K_2K_3)^2,$ we take $(\vv, s)=(\vv_q, s_q(T))$ for a fixed $\vv_q\in \Theta_q$ and $s_q(T):= s_{\vv_q}(\ll_q)\land (T-r_0)$. By the definition of $\Phi_q$, there exists two positive constants $c_1(q)$ and $c_2(q)$ such that
\beg{equation*}\beg{split} \Phi_q(T,\xi,\eta)&\le c_1(q)\Big(1+ \|\xi-\eta\|_\infty^2+ \ff{|\xi(0)-\eta(0)|^2}{c_2(q)\land (T-r_0)}\Big)\\
&\le \ff{c_1(q)(1+c_2(q))}{c_2(q)} \Big(1+\|\xi-\eta\|_\infty^2+\ff{|\xi(0)-\eta(0)|^2}{T-r_0}\Big),\ \ T>r_0, \xi,\eta\in\C.\end{split}\end{equation*}So, for any $p>(1+K_2K_3)^2$ and any $q\in ((1+K_2K_3)^2, p]$, by the first assertion and using the Jensen inequality, we obtain
\beg{equation*}\beg{split} P_T f(\eta) &\le \big(P_Tf^q)^{1/q}(\xi)\exp\bigg[\ff{c_1(q)(1+c_2(q))}{c_2(q)} \Big(1+\|\xi-\eta\|_\infty^2+\ff{|\xi(0)-\eta(0)|^2}{T-r_0}\Big)\bigg]\\
&\le (P_Tf^p)^{1/p}(\xi) \exp\bigg[\ff{c_1(q)(1+c_2(q))}{c_2(q)} \Big(1+\|\xi-\eta\|_\infty^2+\ff{|\xi(0)-\eta(0)|^2}{T-r_0}\Big)\bigg].\end{split}\end{equation*} Therefore, the second assertion holds for
$$C(p)=\inf_{q\in ((1+K_2K_3)^2, p]} \ff{c_1(q)(1+c_2(q))}{c_2(q)} $$ which is decreasing in $p$.

(b) To prove the first assertion, let us fix $\vv\in \Theta_p$ and $t_0\in (0, s_\vv(\ll_p)\land (T-r_0)].$ We shall make use of the coupling constructed in Section 2 for  $\theta=2(1-\vv).$ Since $t_0\le T-r_0$ and $X(t)=Y(t)$ for $t\ge t_0$, we have $X_T=Y_T$ and
\beq\label{4.2} P_T f(\eta)= \E[R_T f(Y_T)]=\E[R_T f(X_T)] \le (P_T f^p(\xi))^{1/p} (\E R_T^{p/(p-1)})^{(p-1)/p}.\end{equation}  By the definition of $R_T$ and $\Q_T$, we have
\beg{equation*}\beg{split} & \E R_T^{p/(p-1)} = \E_{\Q_T} R_T^{1/(p-1)} =\E_{\Q_T} \exp\bigg[\ff 1 {p-1} \int_0^T \<\phi_t,\d \tt B(t)\>
+\ff 1 {2(p-1)}\int_0^T |\phi_t|^2\d t\bigg]\\
&= \E_{\Q_T} \exp\bigg[\ff 1 {p-1}\int_0^T\int_0^T \<\phi_t,\d \tt B(t)\>-\ff{\ss p +1}{2(p-1)^2}\int_0^T|\phi_t|^2\d t +\ff{p+\ss p}{2(p-1)^2}\int_0^T|\phi_t|^2\d t\bigg]\\
&\le \bigg(\E_{\Q_T} \exp\bigg[\ff{\ss p+1}{p-1} \int_0^T \<\phi_t,\d \tt B(t)\>-\ff{(\ss p+1)^2}{2(p-1)^2}\int_0^T |\phi_t|^2\d t\bigg]
\bigg)^{1/(1+\ss p)}\\&\quad\times \bigg(\E_{\Q_T} \exp\bigg[\ff{(\ss p+1)(p+\ss p)}{2(p-1)^2\ss p} \int_0^T |\phi_t|^2 \d t\bigg]\bigg)^{\ss p/(\ss p+1)}\\
&= \bigg(\E_{\Q_T}\exp\bigg[\ll_p\int_0^T |\phi_t|^2\d t\bigg]\bigg)^{\ss p/(\ss p+1)}.\end{split}\end{equation*} Combining this with (\ref{4.2}), we obtain
$$P_T f(\eta)\le (P_T f^p)^{1/p}(\xi) \bigg(\E_{\Q_T}\exp\bigg[\ll_p\int_0^T |\phi_t|^2\d t\bigg]\bigg)^{(\ss p-1)/\ss p}.$$ Therefore, to prove the first assertion, it suffices to show that
\beq\label{4.3}\beg{split}  &\E_{\Q_T} \exp\bigg[\ll_p \int_0^T |\phi_t|^2\d t\bigg]\\
 &\le \exp\bigg[\ff\vv {2(1+\vv)} +\ff{\ll_p(1+\vv)^2K_3^2K_4|\xi(0)-\eta(0)|^2}{2\vv (1-\vv)^2(1+2\vv)(1-\e^{-K_4s})}\\ &\qquad\qquad + \ff{16K_2^2s^2W_{\vv}(\ll_p)}{1-4K_1K_2s}+ \big(K_1^2r_0\ll_p+2sW_{\vv}(\ll_p)\big)\|\xi-\eta\|_\infty^2\bigg].\end{split}\end{equation}Since $X(t)=Y(t)$ for $t\ge t_0$, it is easy to see from the definition of $\phi_t$ and $(A1)$, $(A3)$ that

$$\int_0^T|\phi_t|^2\d t\le  \int_0^{t_0} \bigg\{\ff{K_1^2(1+\vv)}\vv\|X_t-Y_t\|_\infty^2 +\ff{K_2^2(1+\vv)}{\gg(t)^2} |X(t)-Y(t)|^2 \bigg\}\d t +K_1^2r_0\|X_{t_0}-Y_{t_0}\|_\infty^2.$$By this and the H\"older inequality, we obtain
\beq\label{4.4} \beg{split}  &\E_{\Q_T} \exp\bigg[\ll_p \int_0^T |\phi_t|^2\d t\big]]\\
 \le & \bigg(\E_{\Q_T} \exp\bigg[\ll_p K_3^2(1+\vv)^2\int_0^{t_0}\ff{|X(t)-Y(t)|^2}{\gg(t)^2}\d t\bigg]\bigg)^{1/(1+\vv)}\\
&\times \bigg(\E_{\Q_T} \exp\bigg[\ff{2K_1^2(1+\vv)^2\ll_p}{\vv^2} \int_0^{t_0} \|X_t-Y_t\|_\infty^2\d t\bigg]\bigg)^{\vv/(2+2\vv)}\\
&\times\bigg(\E_{\Q_T}
\exp\bigg[\ff{2K_1^2r_0(1+\vv)\ll_p}\vv \|X_{t_0}-Y_{t_0}\|_\infty^2\bigg]\bigg)^{\vv/(2+2\vv)}.\end{split}\end{equation}  Since $\vv\in \Theta_p$ implies that
$$\ll_p K_3^2 (1+\vv)^2 \le \ff{(1-\vv)^4}{2(1+\vv)K_2^2},$$ it follows from Lemma \ref{L4.1} below that
\beq\label{4.5} \beg{split} &\E_{\Q_T} \exp\bigg[\ll_p K_3^2(1+\vv)^2\int_0^{t_0}\ff{|X(t)-Y(t)|^2}{\gg(t)^2}\d t\bigg]\le \exp\bigg[\ff{\ll_pK_3^2(1+\vv)^3|\xi(0)-\eta(0)|^2}{(1+2\vv)(1-\vv)^2\gg(0)}\bigg]\\
&\qquad\qquad\times\bigg(\E_{\Q_T} \exp\bigg[\ff{K_1^2K_2^2K_3^2\ll_p(1+\vv)^3}{8\vv^2(1-\vv)^3}\int_0^{t_0} \|X_t-Y_t\|_\infty^2\d t\bigg]\bigg)^{\vv/(1+2\vv)}.\end{split}\end{equation}Moreover, according to Lemma \ref{L4.2} below,
\beg{equation*}\beg{split} &\E_{\Q_T}
\exp\bigg[\ff{2K_1^2r_0(1+\vv)\ll_p}\vv \|X_{t_0}-Y_{t_0}\|_\infty^2\bigg]\le \exp\bigg[1+\ff{2K_1^2r_0(1+\vv)\ll_p}\vv \|\xi-\eta\|_\infty^2\bigg]\\
&\qquad \times \bigg(\E_{\Q_T}\exp\bigg[\ff{8K_1^3K_2r_0(1+\vv)\ll_p(4K_2K_1r_0(1+\vv)\ll_p +\vv)}{\vv^2} \int_0^{t_0} \|X_t-Y_t\|_\infty^2\d t\bigg]\bigg)^{1/2}.\end{split}\end{equation*} Substituting this and (\ref{4.5}) into (\ref{4.4}), and using the definition of $W_\vv(\ll_p)$, we conclude that
\beq\label{4.6}\beg{split}  &\E_{\Q_T} \exp\bigg[\ll_p \int_0^T|\phi_t|^2\d t\bigg]\le\E_{\Q_T}\exp\bigg[ W_\vv(\ll_p)\int_0^{t_0}\|X_t-Y_t\|_\infty^2\d t\bigg]\\
 & \qquad\times \exp\bigg[\ff{\ll_pK_3^2(1+\vv)^2|\xi(0)-\eta(0|^2}{(1+2\vv)(1-\vv)^2\gg(0)} +\ff\vv{2(1+\vv)} +K_1^2r_0\ll_p\|\xi-\eta\|_\infty^2\bigg].\end{split}\end{equation}Since $t_0\le s_\vv(\ll_p),$ we have
$$W_\vv(\ll_p)\le \ff{1-4K_1K_2t_0}{8K_2^2t_0^2}.$$ So, combining (\ref{4.6}) with Lemma \ref{L4.3} below and noting that for $\theta=2(1-\vv)$ one has
$$\gg(0)= \ff{2\vv}{K_4}(1-\e^{-K_4t_0}),$$  we prove (\ref{4.3}).\end{proof}

\beg{lem}\label{L4.1} For any positive $\ll\le \ff{(1-\vv)^4}{2K_2^2(1+\vv)} $ and $s\in [0,t_0]$,
\beg{equation*}\beg{split} &\E_{\Q_T} \exp\bigg[\ll\int_0^s\ff{|X(t)-Y(t)|^2}{\gg(t)^2}\d t\bigg]\\
&\le \exp\bigg[\ff{\ll(1+\vv)|\xi(0)-\eta(0)|^2}{(1+2\vv)(1-\vv)^2\gg(0)}\bigg] \bigg(\E_{\Q_T} \exp\bigg[\ff{K_1^2K_2^2(1+\vv)\ll}{8\vv^2(1-\vv)^3}\int_0^s\|X_t-Y_t\|_\infty^2\d t\bigg]\bigg)^{\vv/(1+2\vv)}.\end{split}\end{equation*}\end{lem}

\beg{proof} Since $\theta=2(1-\vv)$ and
$$\ff{K_1K_2}{\gg(t)}\|X_t-Y_t\|_\infty|X(t)-Y(t)|\le \ff{K_1^2K_2^2}{4\theta\vv} \|X_t-Y_t\|_\infty^2 +\theta\vv\ff{|X(t)-Y(t)|^2}{\gg(t)^2},$$ it follows from (\ref{*7}) that
$$0\le \tt M(s) +\ff{|\xi(0)-\eta(0)|}{\gg(0)} +\int_0^s \bigg\{\ff{K_1^2K_2^2\|X_t-Y_t\|_\infty^2}{8\vv(1-\vv)} - \ff{2(1-\vv)^2 |X(t)-Y(t)|^2}{\gg(t)^2} \bigg\}\d t.$$Combining this with $(A2)$ and the fact that
\beq\label{FF}\E_{\Q_T} \e^{N(s)+L}\le \big(\E_{\Q_T}\e^{2\<N\>(s)+2L}\big)^{1/2}\end{equation} holds for a $\Q_T$-martingale $N$ and a random variable $L$, we obtain
\beg{equation*}\beg{split} &\E_{\Q_T} \exp\bigg[\ll\int_0^s\ff{|X(t)-Y(t)|^2}{\gg(t)^2}\d t -\ff{\ll |\xi(0)-\eta(0)|^2}{2\gg(0)(1-\vv)^2}\bigg]\\
&\le \E_{\Q_T} \exp\bigg[\ff{\ll}{2(1-\vv)^2}\tt M(s)+\ff{K_1^2K_2^2\ll}{16 \vv (1-\vv)^3}\int_0^s\|X_t-Y_t\|_\infty^2\d t\bigg]\\
&\le \bigg(\E_{\Q_T}\exp\bigg[\ff{2K_2^2\ll^2}{(1-\vv)^4}\int_0^s \ff{|X(t)-Y(t)|^2}{\gg(t)^2}\d t+ \ff {K_1^2K_2^2\ll}{8\vv(1-\vv)^3}\int_0^s\|X_t-Y_t\|_\infty^2\d t\bigg]\bigg)^{1/2}\\
&\le \bigg(\E_{\Q_T}\exp\bigg[\ff{2K_2^2(1+\vv)\ll^2}{(1-\vv)^4} \int_0^s \ff{|X(t)-Y(t)|^2}{\gg(t)^2}\d t\bigg]\bigg)^{1/(2+2\vv)} \\
&\qquad\times \bigg(\E_{\Q_T}\exp\bigg[\ff{K_1^2K_2^2(1+\vv)\ll}{8\vv^2(1-\vv)^3}\int_0^s\|X_t-Y_t\|_\infty^2\d t\bigg]\bigg)^{\vv/(2+2\vv)}.\end{split}\end{equation*} Since
$$\ff{2K_2^2(1+\vv)\ll^2}{(1-\vv)^4}\le\ll$$ and up to an approximation argument as in \cite[Proof of Lemma 2.2]{W10} we may assume that
$$\E_{\Q_T}\exp\bigg[\ll\int_0^s\ff{|X(t)-Y(t)|^2}{\gg(t)^2}\d t\bigg]<\infty,$$   this implies the desired inequality.
\end{proof}

\beg{lem}\label{L4.2} For any $\ll>0$ and $s\in [0,t_0],$
$$\E_{\Q_T} \e^{\ll\|X_s-Y_s\|_\infty^2}\le \e^{1+\ll\|\xi-\eta\|_\infty^2}\bigg(\E_{\Q_T} \exp\bigg[4\ll K_2(2\ll K_2+K_1)\int_0^s\|X_t-Y_t\|_\infty^2\d t\bigg]\bigg)^{1/2}.$$\end{lem}

\beg{proof} Let
$$N(t)= 2\int_0^t \<X(r)-Y(r), (\si(r,X(r))-\si(r,Y(r)))\d\tt B(r)\>,\ \ r\le s,$$ which is a $\Q_T$-martingale. By (\ref{*4}) and noting that $K_4\le \ff 2 {\gg(r)}$, we obtain
\beg{equation*}\beg{split} \|X_t-Y_t\|_\infty^2 &\le \Big\{\sup_{r\in [0,t]}|X(r)-Y(r)|^2\Big\}\lor \|\xi-\eta\|_\infty^2\\
&\le \|\xi-\eta\|_\infty^2 +\sup_{r\in [0,t]} \bigg\{N(r)+ 2K_1K_2\int_0^r\|X_u-Y_u\|_\infty^2\d u\bigg\}.\end{split}\end{equation*}Combining this with (\ref{FF}) and noting that the Doob inequality implies
\beg{equation*}\beg{split} &\E_{\Q_T}\sup_{r\in [0,t]}\e^{M(r)}=\lim_{p\to\infty}\E_{\Q_T}\Big(\sup_{r\in [0,t]}\e^{M(r)/p}\Big)^p\\
&\le \lim_{p\to\infty}\Big(\ff{p}{p-1}\Big)^p \E_{\Q_T}\big(\e^{M(t)/p}\big)^p = \e \,  \E_{\Q_T}\e^{M(t)}\end{split}\end{equation*} for a $\Q_T$-submartingale $M(r)$, we arrive at
\beg{equation*}\beg{split} &\E_{\Q_T} \e^{\ll\|X_s-Y_s\|_\infty^2- \ll \|\xi-\eta\|_\infty^2} \le \E_{\Q_T} \sup_{t\in [0,s]}
\exp\bigg[\ll N(t)+2\ll K_1K_2\int_0^t\|X_r-Y_r\|_\infty^2\d r\bigg]\\
&\le \e\,  \E_{\Q_T} \exp\bigg[\ll N(s)+2\ll K_1K_2\int_0^s\|X_t-Y_t\|_\infty^2\d t\bigg]\\
 &\le \e\, \bigg(\E_{\Q_T}\exp\bigg[2\ll^2\<N\>(s)+4\ll K_1K_2\int_0^s\|X_t-Y_t\|_\infty^2\d t\bigg]\bigg)^{1/2}\\
 &\le \e\,  \bigg(\E_{\Q_T}\exp\bigg[(8K_2^2\ll^2 +4\ll K_1K_2)\int_0^s\|X_t-Y_t\|_\infty^2\d t\bigg]\bigg)^{1/2}.
\end{split}\end{equation*}\end{proof}

\beg{lem}\label{L4.3} For any $s\in (0,t_0]$ and  positive $\ll\le \ff{1-4K_1K_2s}{8K_2^2s^2},$
$$\E_{\Q_T} \exp\bigg[\ll\int_0^s\|X_t-Y_t\|_\infty^2\d t\bigg] \le\exp\bigg[\ff{16K_2^2s^2\ll}{1-4K_1K_2s}+ 2s\ll\|\xi-\eta\|_\infty^2\bigg].$$\end{lem}

\beg{proof} Let
$$\ll_0= \ff{1-4K_1K_2s}{8K_2^2s^2},$$ which is positive since $s\in (0, s_\vv(\ll_p)].$ It is easy to see that
$$4K_2s\ll_0(2K_2s\ll_0+K_1)=\ll_0.$$ So, it follows from the Jensen inequality and Lemma \ref{L4.2} that
\beg{equation*}\beg{split} &\E_{\Q_T} \exp\bigg[\ll_0 \int_0^s\|X_t-Y_t\|_\infty^2\d t\bigg] \le \ff 1 s \int_0^s \E_{\Q_T} \e^{\ll_0 s\|X_t-Y_t\|_\infty^2}\d s\\
&\le \e^{1+\ll_0 s \|\xi-\eta\|_\infty^2}\bigg(\E_{\Q_T} \exp\bigg[4\ll_0K_2s(2\ll_0K_2s+K_1)\int_0^s\|X_t-Y_t\|_\infty^2\d t\bigg]\bigg)^{1/2}\\
&=\e^{1+\ll_0 s \|\xi-\eta\|_\infty^2}\bigg(\E_{\Q_T} \exp\bigg[\ll_0\int_0^s\|X_t-Y_t\|_\infty^2\d t\bigg]\bigg)^{1/2}.\end{split}\end{equation*}Up to an approximation argument as in \cite[Proof of Lemma 2.2]{W10}, we may assume that
$$\E_{\Q_T} \exp\bigg[\ll_0 \int_0^s\|X_t-Y_t\|_\infty^2\d t\bigg] <\infty,$$ so that this implies
$$\E_{\Q_T} \exp\bigg[\ll_0 \int_0^s\|X_t-Y_t\|_\infty^2\d t\bigg] \le \e^{2+2\ll_0 s \|\xi-\eta\|_\infty^2}.$$ Therefore, by the Jensen inequality, for any $\ll\in [0,\ll_0]$
\beg{equation*}\beg{split} \E_{\Q_T} \exp\bigg[\ll_0 \int_0^s\|X_t-Y_t\|_\infty^2\d t\bigg]&\le \bigg(\E_{\Q_T} \exp\bigg[\ll_0 \int_0^s\|X_t-Y_t\|_\infty^2\d t\bigg] \bigg)^{\ll/\ll_0}\\
&\le
\exp\Big[\ff{2\ll}{\ll_0} +2\ll s \|\xi-\eta\|_\infty^2\Big].\end{split}\end{equation*}
\end{proof}

\end{document}